\documentclass[reqno]{amsart}

\usepackage{amsmath,amsfonts,amssymb,verbatim}
\usepackage{epsfig}
\usepackage{graphicx}

\def\homega{\widehat{\omega}}
\def\tdelta{\tilde{\delta}}

\def\e{\epsilon}

\def\N{{\mathbb N}}

\def\R{{\mathbb R}}
\def\Z{{\mathbb Z}}

\def\cK{{\mathcal K}}

\def\1{{\bf 1}}

\def\eqnn{\begin{eqnarray*}}
\def\eeqnn{\end{eqnarray*}}
\def\eqn{\begin{eqnarray}}
\def\eeqn{\end{eqnarray}}

\def\prf{\begin{proof}}
\def\endprf{\end{proof}}

\theoremstyle{plain}
\newtheorem{theorem}{Theorem}[section]

\newtheorem{lemma}[theorem]{Lemma}
\newtheorem{corollary}[theorem]{Corollary}

\numberwithin{equation}{section}

\begin{document}

\parskip=8pt

\title[Blowup rate for Euler equation]
{A lower bound
on blowup rates for the 3D incompressible Euler equation
and a single exponential Beale-Kato-Majda type estimate}

\author[T. Chen]{Thomas Chen}
\address{T. Chen,  
Department of Mathematics, University of Texas at Austin.}
\email{tc@math.utexas.edu}

\author[N. Pavlovi\'{c}]{Nata\v{s}a Pavlovi\'{c}}
\address{N. Pavlovi\'{c},  
Department of Mathematics, University of Texas at Austin.}
\email{natasa@math.utexas.edu}

%\date{\mscrptdate}

\begin{abstract}
We prove a Beale-Kato-Majda type criterion for the loss of regularity for
solutions of the incompressible 
Euler equations
in $H^{s}(\R^3)$, for $s>\frac52$.  
Instead of double exponential estimates of Beale-Kato-Majda type,
we obtain a single exponential bound on $\|u(t)\|_{H^s}$
involving the length parameter introduced by P. Constantin in \cite{co1}. 
In particular, we derive lower bounds on the blowup rate of such solutions. 
\end{abstract}

\maketitle

\section{Introduction}

In this paper, we revisit the Beale-Kato-Majda criterion for the breakdown of 
smooth solutions to the $3D$ Euler equations. 

More precisely, we consider the 
incompressible Euler equations
\begin{align} 
   & \partial_t u + (u\cdot\nabla )u  + \nabla p \, = \, 0 \label{euler} \\
   & \nabla\cdot u \, = \, 0 \label{euler-divfree}\\ 
   &u(x,0) \, = \, u_0 \label{euler-id}
\end{align} 
for an unknown velocity vector $u(x,t) = (u_i(x,t))_{1 \leq i \leq 3}  
\in {\mathbb R}^{3}$ and pressure $p = p(x,t) \in {\Bbb R}$, for position $x \in {\mathbb R}^{3}$ and 
time $t \in [0,\infty)$.

Existence and uniqueness of local in time solutions to \eqref{euler} -- \eqref{euler-id} 
in the space
\eqn\label{lwp-class}
   C([0,T], H^s) \cap C^1([0,T];H^{s-1}) \,,
\eeqn
has long been known for $s>\frac52$,
see for instance  \cite{ka}.  
However, it is an open problem to determine whether such solutions can 
lose their regularity in finite time. 
An important result that addresses the question of a possible loss of regularity of solutions to Euler equations \eqref{euler} -- \eqref{euler-id}
is the criterion formulated by Beale-Kato-Majda \cite{bekama} in terms of
the $L^{\infty}$ norm of the vorticity
$\omega = \nabla \wedge u$.  
More precisely, Beale-Kato-Majda in \cite{bekama} proved the following theorem: 

\begin{theorem} \label{thm-bkm} 
Let $u$ be a solution to \eqref{euler} --  \eqref{euler-id} in the class \eqref{lwp-class}
for $s\geq3$ integer. Suppose that there exists a time 
$T^*$ such that the solution cannot be continued in the class \eqref{lwp-class} to $T=T^*$. 
If $T^*$ is the first such time, then 
\eqn 
   \int_0^{T^*} \|\omega(\cdot,t)\|_{L^{\infty}} \, dt = \infty.
\eeqn
\end{theorem} 

The theorem is proved with a contradiction argument. Under the assumption 
$$
   \int_0^{T^*} \|\omega(\cdot,t)\|_{L^{\infty}} \, dt < \infty \,,
$$ 
the authors of \cite{bekama} show that $\|u(\cdot,t)\|_{H^s} \leq C_0$
for all $t < T^*$, contradicting 
the hypothesis that $T^*$ is the first time such that the solution cannot be continued to $T=T^*$. 
In particular, Beale-Kato-Majda obtain a double exponential bound for $\|u(\cdot,t)\|_{H^s}$, 
which follows from the following estimates: 

\begin{enumerate} 

\item[{\bf Step 1}] An energy-type bound on $\|u\|_{H^s}$ in terms of $\|Du\|_{L^{\infty}}$, where 
$Du=[\partial_i u_j]_{ij}$ is a $3\times 3$-matrix valued function. More specifically, 
one applies the operator $D^{\alpha}$ to equations \eqref{euler}-\eqref{euler-divfree}, 
where $\alpha$ is an integer-valued multi-index  with $|\alpha| \leq s$ and uses a certain 
commutator estimate to derive
\eqn  \label{eq-intro-Hbound} 
   \frac{d}{dt} \| u(\cdot, t)\|_{H^s}^2 \leq 2C \|Du\|_{L^{\infty}} \|u(\cdot, t)\|_{H^s}^2,
\eeqn
which via Gronwall's inequality gives the bound: 
\eqn \label{eq-intro-Hbound-int} 
   \|u(\cdot,t)\|_{H^s} \, \leq \, 
   \|u_0\|_{H^s} \, \exp \left( \, C \int_0^t \|Du(\cdot, \tau)\|_{L^{\infty}} \, d\tau \, \right).
\eeqn
$\;$

\item[{\bf Step 2}] An estimate on $\|Du(\cdot, t)\|_{L^{\infty}}$ based on the quantities
$\|\omega(\cdot, t)\|_{L^{\infty}}$, $\|\omega(\cdot, t)\|_{L^2}$, 
and $\log^{+} \|u(\cdot, t)\|_{H^3}$, given by
\eqn \label{eq-intro-log} 
   \|Du(\cdot,t)\|_{L^{\infty}}\,  \leq \, 
   C \, \left\{ \, 1 \, + \,  \big( \, 1 \, + \, \log^{+}\|u(\cdot,t)\|_{H^3} \, \big) \, \|\omega(\cdot, t)\|_{L^{\infty}} 
   \, + \,  \|\omega(\cdot, t)\|_{L^2} \, \right\}, 
\eeqn
where $C$ is a universal constant. 
\\

\item[{\bf Step 3}] The bound on  $\|\omega(\cdot, t)\|_{L^2}$ in terms of 
$\|\omega(\cdot, t)\|_{L^{\infty}}$ given by
$$
   \frac{d}{dt} \|\omega(\cdot, t) \|_{L^2}^2 
   \, \leq \, 2 \, \widetilde C
   \, \|\omega(\cdot, t)\|_{L^{\infty}} \, \|\omega(\cdot, t) \|_{L^2}^2 \,,
$$ 
which follows from taking the $L^2({\mathbb{R}}^3)$-inner product of
$\omega$ with  the equation for vorticity. 
Then, Gronwall's inequality yields
\eqn \label{eq-intro-vor} 
   \|\omega(\cdot,t)\|_{L^2} \leq \|\omega(\cdot, 0)\|_{L^2} \, 
   \exp \left(\widetilde C  \int_0^t \|\omega(\,\cdot\, , \tau)\|_{L^{\infty}} \, d\tau \right). 
\eeqn
Consequently, one obtains the double exponential bound
\eqn\label{eq-BKMdoubexp-1}
   \|u(\cdot,t)\|_{H^s} \, \leq \, \|u_0\|_{H^s} \, \exp\left( \,
   \exp\Big( \widetilde C \int_0^t \|\omega(\,\cdot\,,\tau)\|_{L^{\infty}} \, d\tau \Big) \, \right) \,.
\eeqn
from combining 
\eqref{eq-intro-Hbound-int}, \eqref{eq-intro-log} and  \eqref{eq-intro-vor}.
\end{enumerate}

It is an open question whether 
\eqref{eq-BKMdoubexp-1} is sharp.
While we do not attempt to answer that question itself in this paper, we obtain 
a single exponential bound on the $H^s$-norm of solution to Euler equations 
\eqref{euler} - \eqref{euler-id} in terms of the quantity
\eqn\label{eq-elldelta-def-1}
    \ell_\delta(t) \, = \, 
    \min\Big\{ \, L \; , \; \Big(\frac{\|\omega(t)\|_{C^\delta}}
    {\|u_0\|_{L^2}}\Big)^{-\frac{2}{2\delta+5}} \, \Big\}, 
\eeqn
where 
\eqn
   \|\omega\|_{C^\delta} \, = \, \sup_{|x-y|<L}\frac{|\omega(x)-\omega(y)|}{|x-y|^\delta}
\eeqn
denotes the $\delta$-Holder seminorm, for 
$L>0$ fixed, and $\delta>0$. More precisely, we prove the following theorem: 

\begin{theorem}\label{thm-main-1}
Let $u$ be a solution to \eqref{euler} - \eqref{euler-id} in the class \eqref{lwp-class},
for $s=\frac52+\delta$. 
Assume that $ \ell_\delta(t) $ is defined as above, and that
\eqn
   \int_0^T \left( \ell_\delta(\tau) \right)^{-\frac52} \, d\tau \, < \, \infty \,.
\eeqn
Then, there exists a finite positive constant $C_\delta=O(\delta^{-1})$ independent
of $u$ and $t$ such that
\eqn 
   \|u(\cdot, t)\|_{H^s}  
   \, \leq \,   \|u_0\|_{H^s} \,
   \exp\Big[ \,  
   C_\delta \, \|u_0\|_{L^2} \, \int_0^t \left( \ell_\delta(\tau) \right)^{-\frac52} \, d\tau \, \Big]  
   \label{eq-main-thm-2-bd-1}
\eeqn
holds for $0 \leq t \leq T$. 
\end{theorem}

The quantity $\ell_{\delta}(t)$ has the dimension of length,  and
was introduced by Constantin in \cite{co1} (see also the work of 
Constantin, Fefferman and Majda \cite{cofema} where
a criterion for loss of regularity in terms of the direction of vorticity was obtained), 
where it was observed that 
\eqn \label{l-cond} 
\int_0^T \left( \ell_\delta(t) \right)^{-\frac{5}{2}} \, dt \, = \, \infty
\eeqn
is a necessary and sufficient condition for blow-up of Euler equations.
In particular, the necessity of the condition follows from the 
inequality obtained in \cite{co1}
\eqn \label{holder-co1} 
\| \omega(\cdot, t) \|_{L^{\infty}}  \leq \|u(\cdot, t)\|_{L^2} \, \left( \ell_\delta(t) \right)^{-\frac{5}{2}}, 
\eeqn 
and Theorem \ref{thm-bkm} of Beale-Kato-Majda.
This is so because Theorem \ref{thm-bkm} implies that if the solution cannot be continued 
to some time $T$, then  $\int_0^T \|\omega(\cdot,t)\|_{L^{\infty}} \, dt = \infty$.
As a consequence of \eqref{holder-co1}, and conservation of energy 
\eqn\label{eq-L2conserv-1}
   \|u(\cdot,t)\|_{L^2} \, = \, \|u_0\|_{L^2} \,,
\eeqn
this in turn implies \eqref{l-cond}. However, by invoking 
the result of Beale-Kato-Majda in this argument, one again obtains a double exponential 
bound on $\|u(\cdot, t)\|_{H^s}$ in terms of $\int_0^T \left( \ell_\delta(t) \right)^{-\frac{5}{2}}dt$. 
We refer to \cite{co2,dehoyu} for recent developments  in this and related areas.

In this paper, we observe that one can actually obtain a single exponential bound on 
the $H^s$-norm of the solution $u(t)$ in terms of $\int_0^T \left( \ell_\delta(t) \right)^{-\frac{5}{2}}dt$, 
as stated in Theorem \ref{thm-main-1}. This is achieved by
avoiding the use of the logarithmic 
inequality \eqref{eq-intro-log} from \cite{bekama}. More precisely, we combine the energy bound 
\eqref{eq-intro-Hbound} with a Calderon-Zygmund type bound on the symmetric 
and antisymmetric parts of $Du$. 

We note that for $s\in\N$, the estimate \eqref{eq-main-thm-2-bd-1}
follows directly from combining Theorem 1 in G. Ponce's paper \cite{po1} with our
Lemma \ref{lm-singintop-1} below, which established the link between 
$\| Du^{+}\|_{L^{\infty}}$ and the length scale $\ell_{\delta}$ (as stated in Corollary 
\ref{cor-Dupm-bds-1}).

As a second main result in this paper, we obtain a lower bound on the blowup rate of solutions
in $H^{\frac52+\delta}$, for an arbitrary, real-valued $\delta>0$. Specifically, we prove: 

\begin{theorem}\label{thm-mainthm-2}
Let $u$ be a solution to \eqref{euler} -- \eqref{euler-id} in the class 
\eqn \label{rate-class}
C([0,T];H^{\frac52 + \delta}) \cap C^1([0,T];H^{\frac32 + \delta}).
\eeqn
Suppose that there exists a time 
$T^*$ such that the solution cannot be continued in the class 
\eqref{rate-class} to $T=T^*$. If $T^*$ is the first such time 
then there exists a finite, positive constant $C(\delta, \|u_0\|_{L^2})$
such that
\eqn\label{eq-blowuprate-1}
   \|u(\cdot, t)\|_{H^{\frac52+\delta}} \, \geq \, C(\delta, \|u_0\|_{L^2}) \, 
   \Big(\frac{1}{T^*-t}\Big)^{1+\frac{2}{5}\delta} \,,
\eeqn
for all $t$ sufficiently close to $T^*$ 
(see the conditions \eqref{eq-smallexp-1} and \eqref{eq-smallexp-2} below, with $t_0 = t$).
\end{theorem}

The proof of Theorem \ref{thm-mainthm-2} can be outlined as follows.
We assume that $u$ is a solution in the class \eqref{rate-class} that cannot be 
continued to $T=T^*$, and that $T^*$ is the first such time. 
Invoking the local in time existence result, 
we derive a lower bound $T_{loc,t_1}>0$ on the time of existence of solutions
to Euler equations in  \eqref{rate-class} for initial data $u(t_1)\in H^{\frac52+\delta}$
at an arbitrary time $t_1 < T^*$.
By definition of $T^*$, we thus have
\eqn \label{eq-intro-times} 
   t_1 \,+\, T_{loc,t_1}  \, <  \, T^* \,.
\eeqn  
Based on an energy bound on the $H^{\frac52 + \delta}$-norm 
of the solution, we obtain in Section \ref{sec-rate} an expression for 
$T_{loc,t_1}$ of the form $\frac{1}{ C \|u(\cdot, t_1)\|_{H^{\frac52+\delta}} }$, which 
together with \eqref{eq-intro-times} implies that 
\eqn \label{eq-intro-trivrate} 
\|u(\cdot, t_1)\|_{H^{\frac52 + \delta}} > \frac{1}{C (T^* - t_1)}\,, 
\eeqn 
for all $t_1 < T^*$. This is an ``a priori" lower bound on the
blowup rate. Subsequently, we improve \eqref{eq-intro-trivrate}  
by a recursion argument in Theorem \ref{thm-mainthm-2} 
for times $t$ close to $T^*$, to yield the stronger bound \eqref{eq-blowuprate-1}.

After completing this work, V. Vicol called to our attention that 
in a recent work, D. Chae proved in \cite{chae} 
(see Theorem 1.1 part (i) of \cite{chae}) that  for 
integer values of $s\in\N$ with $s>1+\frac d2$, and in dimensions $d\geq2$,
\eqn
	\liminf_{t\rightarrow T^*}(T^*-t)\|D^s u(t)\|_{L^2(\R^d)}^{\frac{d+2}{s}}
	\, \geq \, \frac{K}{\|u_0\|_{L^2}^{\frac{d+2}{s}}} 
\eeqn
is a necessary and sufficient condition for blowup at time $T^*$,
where $K=K(d,s)$ is an absolute constant.
In our estimate \eqref{eq-blowuprate-1}, we allow for 
real values of $s=\frac52+\delta$, $\delta>0$, and provide a pointwise lower bound
instead of an infimum limit.

\section{Proof of theorem \ref{thm-main-1}}

%We note the vorticity satisfies the following equation: 
%\begin{align} 
%    \partial_t\omega + (u\cdot\nabla)\omega \, = \, (Du)\omega  \label{vor}\\  
%     \omega(x,t) \, = \, \nabla\wedge u_0 \, \label{vor-id}
%\end{align}
%Here, $Du=[\partial_i u_j]_{ij}$ is a $3\times 3$-matrix valued function.

%Since $u$ is divergence-free, we have that
%\eqn
%    u \, = \, -\Delta^{-1}\nabla\wedge \omega \,,
%\eeqn
%noting that generally, $\nabla\wedge(\nabla\wedge u)=\nabla(\nabla\cdot u)-\Delta u$.

First we recall that the full gradient of velocity $Du$ can be decomposed 
into symmetric and antisymmetric parts,
\eqn
   Du \, = \, Du^+ \, + \, Du^-
\eeqn
where
\eqn
   Du^\pm \, = \, \frac12 \big( \,  Du \pm Du^T \, \big) \,.
\eeqn
$Du^+$ is called the deformation tensor.

In the following lemma, we recall some important properties
of $Du^+$ and $Du^-$. For the convenience of the reader, we give detailed proofs of 
those properties, although they are in part available in the literature, see e.g. \cite{co1}.

\begin{lemma}\label{lm-Dupm-kernel-1}
For both the symmetric and  antisymmetric parts $Du^+$, $Du^-$ of $Du$, 
the $L^2$ bound
\eqn
   \|Du^\pm\|_{L^2} \, \leq \, C  \|\omega\|_{L^2} \,.
\eeqn
holds.

The antisymmetric part $Du^-$ satisfies
\eqn
   Du^- v \, = \, \frac12 \, \omega\wedge v
\eeqn
for any vector $v\in\R^3$.
The vorticity $\omega$ satisfies the identity
\eqn\label{eq-om-Projid-1}
   \omega(x) \, = \, \frac1{4\pi} P.V.\int \sigma(\widehat y)
   \, \omega(x+y) \, \frac{dy}{|y|^3} \,,
\eeqn
("P.V." denotes principal value) 
where
$\sigma(\widehat y)=3\, \widehat y\otimes\widehat y - \1$, with $\widehat y=\frac{y}{|y|}$.
Notably,
\eqn
   \int_{S^2}\sigma(\widehat y) \, d\mu_{S^2}(y) \, = \, 0 \,,
\eeqn
where $d\mu_{S^2}$ denotes the standard measure on the sphere $S^2$.

The matrix components of the symmetric part have the form 
\eqn\label{eq-TijSingIntOp-def-1}
   Du^+_{ij} \, = \,\sum_\ell T_{ij}^\ell(\omega_\ell) \, = \, \sum_\ell \cK_{ij}^\ell*\omega_\ell \,,
\eeqn
where $\omega_\ell$ are the vector components of $\omega$,
and where the integral kernels  $\cK_{ij}^\ell$ have the properties
\eqn\label{eq-cK-prop-1}
   \cK_{ij}^\ell(y)&=&\sigma_{ij}^\ell(\widehat y) \,  |y|^{-3} 
   \\
   \label{eq-cK-prop-2}
   \|\sigma_{ij}^\ell  \|_{C^1(S^2)} &\leq&C 
   \\
   \label{eq-cK-prop-3}
   \int_{S^2}\sigma_{ij}^\ell(\widehat y) \, d\mu_{S^2}(y)&= & 0 \,.
\eeqn 
Thus in particular, $T_{ij}^\ell$ is a Calderon-Zygmund operator, for every $i,j,\ell\in\{1,2,3\}$.
\end{lemma}

\prf
An explicit calculation shows that the Fourier transform of $Du$ as a function of $\homega$ 
is given by
\eqn
   \widehat{Du}(\xi) 
   \, = \, -[(\partial_i(\Delta^{-1}\nabla\wedge\omega)_j )^{\widehat{\;\;}}(\xi)]_{i,j} 
   \, = \, \widehat G(\xi) \, + \, \widehat H(\xi)
\eeqn
where
\eqn
   \widehat G(\xi)\, := \, \frac{1}{2|\xi|^2}\left[ 
   \begin{array}{ccc}
   \xi_1\xi_2\homega_3-\xi_1\xi_3\homega_2  & -\xi_2\xi_3\homega_2 & \xi_2\xi_3\homega_3 \\
   \xi_1\xi_3\homega_1 & \xi_2\xi_3\homega_1-\xi_1\xi_2\homega_3  & -\xi_1\xi_3\homega_3 \\
   -\xi_1\xi_2\homega_1 & \xi_1\xi_2\homega_2& \xi_1\xi_3\homega_2-\xi_2\xi_3\homega_1
   \end{array}\right]
\eeqn
and
\eqn
   \widehat H(\xi)\, := \, \frac{1}{2|\xi|^2}\left[ 
   \begin{array}{ccc}
   0  & \xi_2^2\homega_3 & -\xi_3^2\homega_2 \\
   -\xi_1^2\homega_3 &0  & \xi_3^2\homega_1 \\
   \xi_1^2\homega_2 &- \xi_2^2\homega_1&0
   \end{array}\right] \,,
\eeqn
using the notation $\homega_j\equiv\homega_j(\xi)$ for brevity.

Clearly, every component of $G$ is given by a sum of Fourier multiplication
operators with symbols
of the form $\frac{\xi_i\xi_j}{|\xi|^2}$, $i\neq j$, applied to a component of $\omega$. 
For instance,
\eqn
	G_{21}(x) \, = \, const. \; P.V. \int \, \widehat y_1 \widehat y_3 \, \omega_1(x+y) \,
	\frac{dy}{|y|^3}  
\eeqn
corresponds to the component $G_{21}$.
It is easy to see that every component $G_{ij}$ is a sum of Calderon-Zygmund operators
applied to components of $\omega$,
with kernel satisfying the asserted properties \eqref{eq-cK-prop-1}
$\sim$ \eqref{eq-cK-prop-3}.
The same is true for the
symmetric part, $G^+=\frac12(G+G^T)$.

The symmetric part of $\widehat H(\xi)$ is given by
\eqn
   \widehat H^+(\xi) \, = \, \frac{1}{2|\xi|^2}\left[ 
   \begin{array}{ccc}
   0  & (\xi_2^2-\xi_1^2)\homega_3 & (\xi_1^2-\xi_3^2)\homega_2 \\
   (\xi_2^2-\xi_1^2)\homega_3 &0  &( \xi_3^2-\xi_2^2)\homega_1 \\
   (\xi_1^2-\xi_3^2)\homega_2 &(\xi_3^2- \xi_2^2)\homega_1&0
   \end{array}\right]
\eeqn
so that each component defines a Fourier multiplication operator with symbol
of the form $\frac{\xi_i^2-\xi_j^2}{|\xi|^2}$, $i\neq j$, 
acting on a component of $\omega$ (with associated 
kernel of the form $\frac{x_i^2-x_j^2}{|x|^{n+2}}$). 
That is, for instance,
\eqn
	H_{12}^+(x) \, = \, const \; P.V. \int (\widehat y_2^2-\widehat y_1^2) \,
	\omega_3(x+y) \frac{dy}{|y|^3} \,.
\eeqn
The properties  \eqref{eq-cK-prop-1} $\sim$  \eqref{eq-cK-prop-3}
follow immediately. 

The Fourier transforms of the integral kernels $\cK_{ij}^\ell$ can be read off from the components
$\widehat G^+_{ij}+\widehat H^+_{ij}$. 
In position space, one finds that $\sigma_{ij}^\ell(\widehat y)$ is obtained 
from a sum of terms proportional to terms of the form $\widehat y_{i_1}\widehat y_{j_1}$
and $(\widehat y_{i_2}^2-\widehat y_{j_2}^2)$. 

For the antisymmetric part $Du^-$, one generally has
$Du^-v=\frac12(\nabla\wedge u)\wedge v$
for any $v\in\R^3$, and from $u=-\Delta^{-1}\nabla\wedge \omega$,
we get $Du^-v=\frac12\omega\wedge v$, using that $\nabla\cdot u=0$.

As a side remark,
we note that while $H^-$ does not by itself exhibit 
the properties  \eqref{eq-cK-prop-1} $\sim$  \eqref{eq-cK-prop-3}, 
it combines with $G^-$ in a suitable manner to yield
the stated properties of $Du^-$, thanks to the condition $\nabla\cdot \omega=0$.
\endprf

Next, Lemma \ref{lm-singintop-1} below provides an upper bound
in terms of the quantity  $\ell_\delta(t)$ on singular integral 
operators applied to $\omega$ of the  type appearing in \eqref{eq-TijSingIntOp-def-1}.
We note that similar bounds were used 
in \cite{co1} and \cite{cofema} for the antisymmetric part $Du^-$.
Here, we observe that they also hold for the symmetric part $Du^+$. 
As shown in \cite{cofema} for $Du^-$, the proof of such a bound invokes standard 
arguments based on
decomposing the singular integral into an inner and outer contribution. 
The inner contribution can be bounded  
based on a certain mean zero property, while the outer part is controlled via integration by parts. 

\begin{lemma}\label{lm-singintop-1}
For $L>0$ fixed, and $\delta>0$, let $\ell_\delta(t) $ be defined as above.
Moreover, let $\omega_\ell$,
$\ell=1,2,3$, denote the components of the vorticity vector $\omega(t)$.
Then, any singular integral operator
\eqn\label{eq-singintop-def-1}
   T\omega_\ell(x) \, = \, \frac{1}{4\pi}P.V. \int \, \sigma_T(\widehat y) \, \omega_\ell(x+y)
   \, \frac{dy}{|y|^3}
\eeqn
with 
\eqn\label{eq-Top-prop-2}
   \int_{S^2}\sigma_T(\widehat y)d\mu_{S^2}(y) \, = \, 0
   \; \; \; \; , \; \; \; \; 
   \|\sigma_T\|_{C^1(S^2)} \, < \, C \;,
\eeqn
satisfies
\eqn
   \|T\omega_\ell\|_{L^\infty} \, \leq \, C(\delta) \, \|u_0\|_{L^2} \,  \ell_\delta(t)^{-\frac52}  
\eeqn
for $\ell\in\{1,2,3\}$, for a constant $C(\delta)=O(\delta^{-1})$ independent of $u$ and $t$.
\end{lemma}

\prf
Let $\chi_1(x)$ be a smooth cutoff function which is identical to 1 on $[0,1]$,
and identically 0 for $x>2$. Moreover, let $\chi_R(x)=\chi_1(x/R)$, and
$\chi_R^c=1-\chi_R$.

We consider
\eqn
   \int_{|y|>\e} \, \sigma_T(\widehat y) \, \omega_\ell(x+y)
   \, \frac{dy}{|y|^3} \, = \, (I) \, + (II)
\eeqn
for $\e>0$ arbitrary, where
\eqn
   (I) \, := \, \int_{|y|>\e} \, \sigma_T(\widehat y) \, \omega_\ell(x+y) \, \chi_{\ell_\delta(t)} (|y|)\,
   \, \frac{dy}{|y|^3}
\eeqn
and
\eqn
   (II) \, := \, \int \, \sigma_T(\widehat y) \, \omega_\ell(x+y) \, \chi_{\ell_\delta(t)} ^c(|y|) \,
   \, \frac{dy}{|y|^3} \,.
\eeqn
From the zero average property \eqref{eq-Top-prop-2}, we find
\eqn\label{eq-termI-est-1}
   |(I)| &=&
   \Big| \, \int_{|y|>\e} \, \sigma_T(\widehat y) \, (\omega_\ell(x+y) -\omega_\ell(x) )
   \, \chi_{\ell_\delta(t)} (|y|)\,
   \, \frac{dy}{|y|^3} \, \Big|
   \nonumber\\
    &\leq&\|\omega_\ell\|_{C^\delta} \,
   \int_{|y|<2\ell_\delta(t) }\frac{dy}{|y|^{3-\delta}} 
   \nonumber\\
   &\leq&\frac C\delta (\ell_\delta(t))^{\delta} \, \|\omega_\ell\|_{C^\delta} 
   \nonumber\\
   &\leq&C \, \delta^{-1} \, \|u_0\|_{L^2} \,  (\ell_\delta(t))^{-\frac52}
\eeqn
since from the definition of $ \ell_\delta(t)$,
\eqn 
    \|\omega_\ell\|_{C^\delta} \, \leq \, \|u_0\|_{L^2} \, (\ell_\delta(t))^{-\delta-\frac52}
\eeqn
follows straightforwardly.
We  can send $\e\searrow0$, since  the estimates are uniform in $\e$.

On the other hand,
\eqn
   (II) &=& \int \, \sigma_T(\widehat y) \, (\partial_{y_i} u_j-\partial_{y_j} u_i)(x+y) 
   \, \chi_{\ell_\delta(t)} ^c(|y|) \,
   \, \frac{dy}{|y|^3} \,.
\eeqn
It suffices to consider one of the terms in the difference,
\eqn
   \lefteqn{
   \Big| \int \, \sigma_T(\widehat y) \, \partial_{y_i} u_j (x+y) 
   \, \chi_{\ell_\delta(t)} ^c(|y|) \,
   \, \frac{dy}{|y|^3} \,  \Big|
   }
   \nonumber\\
   &=&\Big|\int \, dy \, u_j (x+y) 
   \,\partial_{y_i}\Big(  \sigma_T(\widehat y) \, \chi_{\ell_\delta(t)} ^c(|y|) \,
   \, \frac{1}{|y|^3} \, \Big) \, \Big|
   \nonumber\\
   &\leq&
   C \, \|u_j\|_{L^2} \, \Big\| \, \partial_{y_i}\Big(\sigma_T(\widehat y) \, \chi_{\ell_\delta(t)} ^c(|y|) \,
   \, \frac{1}{|y|^3} \, \Big) \, \Big\|_{L^2}
   \nonumber\\
   &\leq&
   C \, \|u_0\|_{L^2} \, (\ell_\delta(t))^{-\frac52}
\eeqn
where to obtain the last line we used the conservation of energy 
\eqref{eq-L2conserv-1} and the following three bounds: 
\begin{enumerate}
\item [(i)]

\eqn 
   \Big\| \,\Big(  \partial_{y_i}\chi_{R} ^c(|y|) \, \Big)  \,
   \, \frac{\sigma_T(\widehat y)}{|y|^3} \, \Big\|_{L^2}^2 
   &\leq&C \,
    \frac1{R^2}\int_{R<|y|<2R} \frac{dy}{|y|^6}
    \nonumber\\
    &\leq&C \, R^{-5} \,,
\eeqn
for $R=\ell_\delta(t)$. 

\item[(ii)]
\eqn 
   \Big\| \, \sigma_T(\widehat y) \, \chi_{R} ^c(|y|) \,   \,
   \, \partial_{y_i}\frac{1}{|y|^3} \, \Big\|_{L^2}^2 
   &\leq&C \,
   \int_{|y|>R} \frac{dy}{|y|^8}
    \nonumber\\
    &\leq&C \, R^{-5} \,.
\eeqn

\item[(iii)]
\eqn 
   \Big\| \,   \chi_{R} ^c(|y|) \,  \frac{1}{|y|^3}   \,
   \, \partial_{y_i}\,  \sigma_T(\widehat y) \, \Big\|_{L^2}^2 
   &\leq&C \,
   \int_{|y|>R}  \frac{1}{|y|^6}   \, \frac{1}{|y|^2}   \, dy 
    \nonumber\\
    &\leq&C \, R^{-5} \,,
\eeqn
where we used that
\eqn
   \Big|\nabla_y\sigma_T(\widehat y)\Big| &=&
   \Big| \, \frac{1}{|y|}
%   (\1-\widehat y\otimes\widehat y) 
   \, (\nabla_z\sigma_T(z_1,z_2,z_3))\Big|_{z=\widehat y}\, \Big|
   \nonumber\\
   &\leq&\frac{1}{|y|} \, \|\sigma_T\|_{C^1(S^2)} 
\eeqn
holds.
\end{enumerate}

Summarizing, we arrive at
\eqn
   \|T\omega_\ell\|_{L^\infty} \, \leq \, C(\delta) \, \|u_0\|_{L^2} \,  \ell_\delta(t)^{-\frac52}  
\eeqn
for $C(\delta)=O(\delta^{-1})$, which is the asserted bound.
\endprf

The form of the singular integral operator that appears in the statement of Lemma \ref{lm-singintop-1}
is suitable for application to $Du^+$ and $Du^-$, as we shall see in the following corollary. 

\begin{corollary}\label{cor-Dupm-bds-1}
There exists a finite, positive constant $C_\delta=O(\frac1\delta)$ independent of $u$ and $t$
such that the estimate 
\eqn
   \|Du^+\|_{L^\infty} +
   \|Du^-\|_{L^\infty} 
   \, \leq \, C_\delta \,  \|u_0\|_{L^2} \, \ell_\delta(t)^{-\frac52}
\eeqn
holds.
\end{corollary}

\prf
According to Lemma \ref{lm-Dupm-kernel-1},
the matrix components of both $Du^+$ and $Du^-$ have the form 
\eqref{eq-singintop-def-1}. 

Accordingly, Lemma \ref{lm-singintop-1}
immediately implies the assertion.
\endprf

Now we are ready to give a proof of Theorem \ref{thm-main-1}, 
which is based on combining an energy estimate for Euler equations with Corollary 
\ref{cor-Dupm-bds-1}. 

\begin{comment} 
Similarly as in \cite{bekama}, we consider 
\eqn
   \frac12 \partial_t\|\omega(t)\|_{L^2}^2 
   & = & \big( \, \omega(t) \, , \, (Du(t)) \, \omega(t) \, \big)_{L^2}
   \nonumber\\
   & = & \big( \, \omega(t) \, , \, (Du^+(t)) \, \omega(t) \, \big)_{L^2}
\eeqn
which implies that 
\eqn
   \frac12 \partial_t\|\omega(t)\|_{L^2}^2 
   & \leq & C \, \|\omega(t)\|_{L^\infty} \, \|\omega(t)\|_{L^2}\|Du^+(t)\|_{L^2}
   \nonumber\\
   &\leq&C \, \|\omega(t)\|_{L^\infty} \, \|\omega(t)\|_{L^2}^2  \,.
\eeqn 
Accordingly,
\eqn\label{eq-omL2-bd-1}
   \|\omega(t)\|_{L^2} \, \leq \, \exp\big( \, \int_0^t\|\omega(s)\|_{L^\infty} \, ds\, \big) 
   \, \|\omega_0\|_{L^2} \,.
\eeqn 
In \cite{bekama}, the strategy of the argument is to bound $\|Du(t)\|_{L^\infty}$
by $\|\omega(t)\|_{L^\infty}$, whereby the authors obtained a logarithmic correction
of the form $\log(\|u(t)\|_{H^3})$. Based on this, the authors applied the Gronwall
inequality to $\|u(t)\|_{H^s}$ with $s\geq3$, thereupon obtaining double exponential
bounds in $t$ (due to Gronwall, and due to the logarithm).

We improve these   estimates by making use of
the special properties of the symmetric part,  $Du^+(t)$, derived above. 
\end{comment}

For $s\geq3$ integer-valued, the energy bound \eqref{eq-intro-Hbound} 
\eqn \label{eq-ineq-sobD}
   \frac12 \partial_t\|u(t)\|_{H^s}^2 \, \leq \,   \| Du(t) \|_{L^\infty} \, \|u(t)\|_{H^s}^2 \,
\eeqn 
was proven in \cite{bekama}. For fractional $s>\frac52$, 
 we recall the definitions of the homogenous and inhomogenous  Besov norms
for $1\leq p,q\leq \infty$,
\eqn
   \|u\|_{\dot B^s_{p,q}}
   \, = \, 
   \Big(\sum_{j\in\Z} 2^{jqs} \|u_j\|_{L^p}^q\,\Big)^{\frac1q} \,,
\eeqn
respectively,
\eqn
   \|u\|_{B^s_{p,q}}
   \, = \, 
   \Big(\| u \|_{L^p}^q \, + \, \|u\|_{\dot B^s_{p,q}}^q \,\Big)^{\frac1q} \,,
\eeqn
where $u_j=P_j u$ is the Paley-Littlewood projection of $u$ of scale $j$.
In analogy to  \eqref{eq-ineq-sobD}, we obtain
the bound on the $B^s_{2,2}$ Besov norm
of $u(t)$ given by
\eqn\label{eq-uBesovsquare-1}
   \frac12 \partial_t\|u(t)\|_{B^s_{2,2}}^2 \, \leq \,   \| Du(t) \|_{L^\infty} \, 
   \|u(t)\|_{B^s_{2,2}}^2 \,,
\eeqn
from a straightforward
application of estimates obtained in \cite{pl}; details\footnote{We note that for integer $s \geq 3$, G. Ponce obtained in \cite{po1} 
the following improvement of \eqref{eq-ineq-sobD},
$$
   \frac12 \partial_t\|u(t)\|_{H^s}^2 \, \leq \,   \| Du^{+}(t) \|_{L^\infty} \, \|u(t)\|_{H^s}^2.
$$ 
Our proof of \eqref{eq-uBesovsquare-1} for fractional  $s>\frac{5}{2}$ does not yield 
the analogous improved bound.}
are given in the Appendix.  
Accordingly,  since the left hand side yields 
\eqn
   \partial_t\|u(t)\|_{B^s_{2,2}}^2 \, = \, 2 \|u(t)\|_{B^s_{2,2}}
   \partial_t\|u(t)\|_{B^s_{2,2}} \,,
\eeqn
we get
\eqn \label{eq-proof1-Hbound} 
   \partial_t\|u(t)\|_{B^s_{2,2}} \, \leq \,   \| Du(t) \|_{L^\infty} \, \|u(t)\|_{B^s_{2,2}} \,.
\eeqn 
However, Corollary \ref{cor-Dupm-bds-1} implies that 
\eqn \label{eq-proof1-cor}
    \| Du(t) \|_{L^\infty} & \leq &  \| Du^+(t) \|_{L^\infty} \, + \,  \| Du^-(t) \|_{L^\infty}
   \nonumber\\
   &\leq&C_\delta \, \|u_0\|_{L^2} \,  \, (\ell_\delta(t))^{-\frac52}   \,.
\eeqn 
Therefore, by combining \eqref{eq-proof1-Hbound} and \eqref{eq-proof1-cor} we obtain 
 \eqn\label{eq-dert-uHalph-1} 
   \partial_t \|u(t)\|_{B^s_{2,2}}
   & \leq & C_\delta \, \|u_0\|_{L^2} \, (\ell_\delta(t))^{-\frac52}   \, \|u(t)\|_{B^s_{2,2}} \,,
   \nonumber
\eeqn
which implies that
\eqn 
   \|u(t)\|_{H^s}  & \sim &
   \|u(t)\|_{B^s_{2,2}}
    \nonumber\\
   & \leq &  \|u_0\|_{B^s_{2,2}} \,\exp\Big[ \,  
   C_\delta \, \|u_0\|_{L^2} \, \int_0^t \ell_\delta(s)^{-\frac52} \, ds \, \Big]
   \nonumber\\
   & \sim &  \|u_0\|_{H^s} \,\exp\Big[ \,  
   C_\delta \, \|u_0\|_{L^2} \, \int_0^t \ell_\delta(s)^{-\frac52} \, ds \, \Big]  \,,
   \nonumber
\eeqn
for $s\geq0$, where we recall 
from \eqref{eq-termI-est-1} that $C_\delta=O(\delta^{-1})$
(see also \cite{po1} for a related bound, but without \eqref{eq-proof1-cor}). 

This completes the proof of Theorem \ref{thm-main-1}.
\qed

\section{Lower bounds on the blowup rate} \label{sec-rate}

In this section, we prove Theorem \ref{thm-mainthm-2}.

Recalling the energy bound \eqref{eq-proof1-Hbound},
\eqn
   \partial_t\|u(t)\|_{B^s_{2,2}} \, \leq \, \|Du(t)\|_{L^\infty} \, \|u(t)\|_{B^s_{2,2}} \,,
\eeqn
we invoke the Sobolev embedding 
\eqn
   \|Du\|_{L^\infty} &\leq & \|\widehat{Du}\|_{L^1} 
   \nonumber\\
   &\leq&\Big(\int d\xi \, \langle \xi\rangle^{-3-2\delta}\Big)^{\frac12}\|Du\|_{H^{\frac32+\delta}}
   \nonumber\\
   &\leq&C_\delta \, \|u\|_{H^{\frac52+\delta}} 
   \nonumber\\
   &\sim& C_\delta \, \|u\|_{B^{\frac52+\delta}_{2,2}} \,,
\eeqn
with $C_\delta=O(\delta^{-{\frac12}})$, to get, for $s=\frac52+\delta$, 
\eqn
   \partial_t \|u(t)\|_{B^s_{2,2}} \, \leq \, C_\delta \, (\|u(t)\|_{B^s_{2,2}})^2 \,.
\eeqn
Straightforward integration implies
\eqn
   -\Big(\frac{1}{\|u(t)\|_{B^s_{2,2}} }-\frac{1}{ \|u(t_0)\|_{B^s_{2,2}} } \Big) 
   \, \leq \, C_\delta(t-t_0)\,.
\eeqn
Hence, 
\eqn
   \|u(t)\|_{H^s} 
   & \sim & 
   \|u(t)\|_{B^s_{2,2}}
   \nonumber\\
   & \leq &
   \frac{ \|u(t_0)\|_{B^s_{2,2}} }{1-(t-t_0) C_\delta  \|u(t_0)\|_{B^s_{2,2}} }
   \nonumber\\
   & \sim &
   \frac{ \|u(t_0)\|_{H^s} }{1-(t-t_0) C_\delta  \|u(t_0)\|_{H^s} } \,,
\eeqn 
where a possible trivial modification of $C_\delta$ is implicit in passing to the last line.
This implies that the solution $u(t)$ is locally well-posed in $H^s$, with $s = \frac52 + \delta$, for 
\eqn
   t_0 \, \leq \, t \, < \, t_0+\frac{1}{C_\delta\|u(t_0)\|_{H^s}} \,.
\eeqn
In particular, this infers that if $T^*$ is the first time beyond which the solution $u$ 
cannot be continued, one necessarily has that
\eqn
   T^*  \, > \,  t_0+\frac{1}{C_\delta\|u(t_0)\|_{H^s}} \,.
\eeqn
This in turn implies an a priori lower bound on the blowup rate given by
\eqn\label{eq-apriori-blwpbd-1}
   \|u(t)\|_{H^s} \, > \, \frac{1}{C_\delta \, (T^*-t)}
\eeqn
for all $0\leq t <T^*$. 

Next, we derive the lower bound on the blowup rate stated in Theorem \ref{thm-mainthm-2}
which is stronger than \eqref{eq-apriori-blwpbd-1}.
To begin with, we note that
\eqn\label{eq-omegaHolder-1}
   \|\omega(t)\|_{C^\delta} &\leq& C_\delta \|\omega(t)\|_{H^{\frac32+\delta}}
   \nonumber\\
   &\leq&C_\delta \|u(t)\|_{H^{\frac52+\delta}}
   \nonumber\\
   &\leq & \frac{C_\delta \|u(t_0)\|_{H^{\frac52+\delta}}}
   {1- (t-t_0) C_\delta  \|u(t_0)\|_{H^{\frac52+\delta}}} \,.
\eeqn
That is, local well-posedness of $u$ in $H^{\frac52+\delta}$ implies $\delta$-Holder
continuity of the vorticity.

The parameter $L$ in the definition \eqref{eq-elldelta-def-1}  of $\ell_\delta(t)$ 
is arbitrary. Thus, in view of \eqref{eq-omegaHolder-1},
we may now let $L\rightarrow\infty$ for convenience. Then,
%Let us now assume that that the $\delta$-Holder
%norm of $\omega(t)$ is sufficiently large that
\eqn\label{eq-ldelta-bd-1}
   \ell_\delta(t)^{-\frac52} &=&
   \Big( \frac{ \|\omega(t)\|_{C^\delta} }{ \|u_0\|_{L^2} } \Big)^{\frac{2}{2\delta+5}\cdot \frac52}
   \nonumber\\
   &\leq&
   \Big(\frac{C_\delta \,  \|u(t)\|_{ H^{\frac52+\delta} } }{\|u_0\|_{L^2} } \Big)^{1-\tdelta}
   \nonumber\\
   &\leq &\Big(\frac{C_\delta}{ \|u_0\|_{L^2} } \Big)^{1-\tdelta}
   \Big(\frac{  \|u(t_0)\|_{H^s}}{1- (t-t_0) C_\delta  \|u(t_0)\|_{H^s}} 
   \Big)^{1-\tdelta} \,,
\eeqn
where
\eqn
   \tdelta \, := \, \frac{2\delta}{5+2\delta} \, {\mbox{ and }} \, s = \frac52+\delta\,.
\eeqn
We note that while the right hand side of \eqref{eq-ldelta-bd-1} 
diverges as $t$ approaches 
\eqn
   t_1 \, := \, t_0 \, + \, \frac{1}{C_\delta\|u(t_0)\|_{H^s}} \,,
\eeqn
the integral 
\eqn
   \int_{t_0}^{t_1}\ell_\delta(t)^{-\frac52}dt
   & \leq & 
   \Big(\frac{C_\delta}{ \|u_0\|_{L^2} } \Big)^{1-\tdelta}
   \int_{t_0}^{t_1}\Big(\frac{  \|u(t_0)\|_{H^s}}{1- (t-t_0) C_\delta  \|u(t_0)\|_{H^s}} 
   \Big)^{1-\tdelta}dt
   \nonumber\\
   & =: & B_0(\delta) \,
\eeqn
{\em converges} for $\delta>0$ ($\Leftrightarrow$ $\widetilde \delta>0$).
This implies that the solution $u(t)$ for $t\in[t_0,t_1)$ can be extended to $t>t_1$.

In particular, we obtain that
\eqn
   \|u(t_1)\|_{H^{\frac52+\delta}} &\leq&
   \|u(t_0)\|_{H^{\frac52+\delta}} \,
   \exp\Big( \, C_\delta \, \|u_0\|_{L^2}\int_{t_0}^{t_1}(\ell_\delta(t))^{-\frac52}dt \, \Big) 
   \nonumber\\
   &\leq& \|u(t_0)\|_{H^{\frac52+\delta}} 
   \, \exp\Big( \, C_\delta \, \|u_0\|_{L^2}B_0(\delta) \, \Big) 
\eeqn
from Theorem \ref{thm-main-1}.

We may now repeat the above estimates with initial data $u(t_1)$ in $H^{\frac52+\delta}$,
thus obtaining a local well-posedness interval $[t_1,t_2]$. 
Accordingly, we may set $t_2$ to be given by
\eqn
   t_2  \, := \, t_1 \, + \, \frac{ 1 }{C_\delta\|u(t_1)\|_{H^s}}   \,.
\eeqn
More generally, we define the discrete times $t_j$ by
\eqn
   t_{j+1} \, := t_j \, + \, \frac{ 1 }{C_\delta\|u(t_j)\|_{H^s}} \,.
\eeqn
We then have 
\eqn
   \|u(t_{j+1})\|_{H^s}  \, \leq \, 
   \exp\Big( \, C_\delta\|u_0\|_{L^2}
   \, B_j(\delta)\, \Big) \, \|u(t_{j})\|_{H^s} \,,
\eeqn
where $B_j(\delta)$ is defined by
\eqn
   \lefteqn{
   C_\delta\|u_0\|_{L^2}B_j(\delta)
   }
   \nonumber\\
   &:=&
   C_\delta\|u_0\|_{L^2}
   \Big(\frac{C_\delta}{ \|u_0\|_{L^2} } \Big)^{1-\tdelta}
   \int_{t_j}^{t_{j+1}}\Big(\frac{  \|u(t_j)\|_{H^s}}{1- (t-t_j) C_\delta  \|u(t_j)\|_{H^s}} 
   \Big)^{1-\tdelta}dt
   \nonumber\\ 
%    &= &C_\delta\|u_0\|_{L^2} \Big(\frac{C_\delta}{ \|u_0\|_{L^2} } \Big)^{1-\tdelta}
%    \int_{t_j}^{t_{j+1}}\Big(\frac{  \|u(t_j)\|_{H^\alpha}}{1- (t-t_j) C_\delta  \|u(t_j)\|_{H^\alpha}} 
%    \Big)^{1-\tdelta}dt
%    \nonumber\\
   &= &\frac{1}{\tdelta} C_\delta^{1-\tdelta} \, 
   \Big(\frac{ \|u_0\|_{L^2} }{ \|u(t_j)\|_{H^s}} \Big)^{\tdelta}
   \nonumber\\
   &=: &b_\delta \,  
   \Big(\frac{ \|u_0\|_{L^2} }{ \|u(t_j)\|_{H^s}} \Big)^{\tdelta} \,.
\eeqn
Letting
\eqn\label{eq-rhoj-def-1}
   \rho_j \, := \, \exp\Big(b_\delta \, 
   \Big(\frac{ \|u_0\|_{L^2} }{ \|u(t_j)\|_{H^s}} \Big)^{\tdelta}\Big) \,,
\eeqn
we have  
\eqn
   \|u(t_{j})\|_{H^s}  \, \leq \, \rho_{j-1}  \, \|u(t_{j-1})\|_{H^s} \,,
\eeqn
and we remark that $(\rho_j)_j$ satisfy the recursive estimates
\eqn
   \rho_j &\geq& \exp\Big(b_\delta \, 
   \Big(\frac{ \|u_0\|_{L^2} }{ \rho_{j-1}\|u(t_{j-1})\|_{H^s}} \Big)^{\tdelta}\Big) 
   \nonumber\\
   &=&
   ( \, \rho_{j-1} \, )^{\rho_{j-1}^{-\tdelta}}
   \nonumber\\
   &=&\exp\Big(\, \rho_{j-1}^{-\tdelta} \, \ln \rho_{j-1} \,\Big)\, .
\eeqn
We note that from its definition, $\rho_j>1$ for all $j$.

We shall now assume that  $T^*>0$ is the first time beyond which
the solution $u(t)$ cannot be continued.
Thus, by choosing $t_0$ close enough to $T^*$, \eqref{eq-apriori-blwpbd-1}
implies that
$\|u(t_0)\|_{H^s}$ can be made sufficiently large that
the following hold:
\begin{enumerate}
\item
The quantity 
\eqn\label{eq-smallexp-1}
   b_\delta \, 
   \Big(\frac{ \|u_0\|_{L^2} }{ \|u(t_0)\|_{H^s}} \Big)^{\tdelta} \, \ll \, 1 \,
\eeqn
is small.
\\
\item
There is a positive, finite constant $\widetilde C$ independent of $j$ such that
\eqn\label{eq-smallexp-2}
   \|u(t_j)\|_{H^s} \, \geq \, \widetilde C \, \|u(t_0)\|_{H^s}
\eeqn 
holds for all $j\in\N$.
Without any loss of generality
(by a redefinition of the constant $b_\delta$ if necessary), 
we can assume that $\widetilde C = 1$.

We note that in principle, there might be strong oscillations 
close to blowup so that \eqref{eq-smallexp-2} is not obviously true.
The fact that  \eqref{eq-smallexp-2} holds follows from
the a priori bound  \eqref{eq-apriori-blwpbd-1}.
\end{enumerate}

Accordingly, \eqref{eq-smallexp-2} with $\widetilde C=1$ implies that 
$\rho_j\leq  \rho_0$ for all $j$.
Then, for any $N\in\N$,
\eqn
   T^* -t_0&\geq&\sum_{j=0}^{N}(t_{j+1}-t_j)
   \nonumber\\
   &=&\frac{1}{C_\delta} \, \Big( \, \frac{1}{\|u(t_0)\|_{H^s}} \, + \cdots +
   \, \frac{1}{\|u(t_N)\|_{H^s}}\, \Big)
   \nonumber\\
   &=&\frac{1}{C_\delta\|u(t_0)\|_{H^s}} 
   \Big( \, 1\,+\, \frac{\|u(t_0)\|_{H^s}}{\|u(t_1)\|_{H^s}} \, + \cdots +
   \, \frac{\|u(t_0)\|_{H^s}}{\|u(t_N)\|_{H^s}}\, \Big)
   \nonumber\\
   &\geq&\frac{1}{C_\delta\|u(t_0)\|_{H^s}} 
   \Big( \, 1\,+\, \frac{1}{\rho_0} \, + \cdots +
   \, \frac{1}{\rho_0\cdots\rho_N}\, \Big)
   \nonumber\\
   &\geq&\frac{1}{C_\delta\|u(t_0)\|_{H^s}} 
   \Big( \, 1\,+\, \frac{1}{\rho_0} \, + \cdots +
   \, \frac{1}{\rho_0^N}\, \Big)
\eeqn
from $\frac{1}{\rho_j} \geq  \frac{1}{\rho_0}$ for all $j$,
and the fact that $\rho_0>1$ since the argument in the exponent \eqref{eq-rhoj-def-1}
is positive.

Then, letting $N\rightarrow\infty$, we obtain
\eqn
   \frac{1}{T^* -t_0}&\leq&C_\delta\|u(t_0)\|_{H^s}\Big( \, 1-\frac1{\rho_0} \, \Big) 
   \nonumber\\
   &=&C_\delta\|u(t_0)\|_{H^s}\Big(\, 1- \exp\Big(-b_\delta \, 
   \Big(\frac{ \|u_0\|_{L^2} }{ \|u(t_0)\|_{H^s}} \Big)^{\tdelta}\Big) \, \Big) \,.
\eeqn
Next, we deduce a lower bound on the blowup rate.

Invoking \eqref{eq-smallexp-1}, we obtain 
\eqn
   \frac{1}{T^* -t_0}&\leq&C_\delta\|u(t_0)\|_{H^s}\Big(\, 1- \exp\Big(-b_\delta \, 
   \Big(\frac{ \|u_0\|_{L^2} }{ \|u(t_0)\|_{H^s}} \Big)^{\tdelta}\Big) \, \Big) 
   \nonumber\\
   &\approx&C_\delta\|u(t_0)\|_{H^s} b_\delta \, 
   \Big(\frac{ \|u_0\|_{L^2} }{ \|u(t_0)\|_{H^s}} \Big)^{\tdelta} 
   \nonumber\\
   &=&C_\delta \, b_\delta \|u_0\|_{L^2}^{\tdelta}
   \|u(t_0)\|_{H^s}^{1-\tdelta}  \, .
\eeqn
This implies a lower bound on the blowup rate of the form
\eqn\label{eq-blowuprate-2}
   \|u(t_0)\|_{H^{\frac52+\delta}} & \geq & C(\delta, \|u_0\|_{L^2}) \, 
   \Big(\frac{1}{T^*-t_0}\Big)^{\frac{1}{1-\tdelta}}
   \nonumber\\
   & = & 
   C(\delta, \|u_0\|_{L^2}) \, 
   \Big(\frac{1}{T^*-t_0}\Big)^{\frac{2\delta+5}{5}} \,,
\eeqn
under the condition that \eqref{eq-smallexp-1}  and \eqref{eq-smallexp-2} hold.

%On the other hand, if the left hand side of \eqref{eq-smallexp-1} is not small
%(far away from the blowup time), then
%we obtain
%\eqn
%    \|u(t)\|_{H^\alpha} \, \geq \, C(\delta, \|u_0\|_{L^2},\|u(t_0)\|_{H^\alpha} ) \, 
%    \frac{1}{T-t} \,.
%\eeqn
%which is slower than \eqref{eq-blowuprate-1}.

This concludes our proof of Theorem \ref{thm-mainthm-2}.
\qed

\appendix

\section{Proof of inequality \eqref{eq-proof1-Hbound} for $s>\frac52$}

In this Appendix, we prove \eqref{eq-proof1-Hbound} which follows from
\eqref{eq-uBesovsquare-1},  
\eqn\label{eq-App-Enbd-1}
   \frac12 \partial_t \|u(t)\|_{B^s_{2,2}}^2 \, \lesssim \, \|D  u(t)  \|_{L^\infty} \, 
   \|u(t)\|_{B^s_{2,2}}^2 \,,
\eeqn
for $s>\frac52$.
We invoke Eq. (26) in the work \cite{pl} of F. Planchon,
which is valid for $s>1+\frac n2$ in $n$ dimensions
(thus, $s>\frac52$ in our case of $n=3$),
for parameter values $p=q=2$ in the notation of that paper. It yields
\eqn
   \frac12\partial_t 2^{2js} \|u_j\|_{L^2}^2 &\lesssim&
   2^{2js} \sum_{k\sim j} \| S_{j+1} D u \|_{L^\infty} \, \|u_k\|_{L^2} \, \|u_j\|_{L^2}
   \nonumber\\
   &&\hspace{1cm} + \, 2^{2js}\sum_{j\lesssim k\sim k'}
    \|u_k\|_{L^2} \,  \|u_{k'}\|_{L^2} \,  \|D u_j\|_{L^\infty} 
\eeqn
where $u_k=P_k u$ is the Paley-Littlewood projection of $u$ at scale $k$, and
$S_j=\sum_{j'\leq j}P_{j'}$ is the Paley-Littlewood projection to scales $\leq j$.

Summing over $j$,
\eqn
   \frac12\partial_t \sum_j 2^{2js} \|u_j\|_{L^2}^2 
   &\lesssim&
   \sup_j \|S_{j+1} D  u \|_{L^\infty} \, \Big(
   \sum_j 2^{2js} \sum_{k\sim j}  \|u_k\|_{L^2} \, \|u_j\|_{L^2}
   \nonumber\\
   &&\hspace{1cm} + \, 
   \sum_j\sum_{k\sim k' \gtrsim j}  2^{2s(j-k)}
   2^{ks} \, \|u_k\|_{L^2} \, 2^{k's} \, \|u_{k'}\|_{L^2} \, \Big)
   \nonumber\\
   &\lesssim&
   \|D  u \|_{L^\infty} \, \Big(
   \sum_j 2^{2js}   \|u_j\|_{L^2}^2
   \nonumber\\
   &&\hspace{1cm} + \, 
   \sum_{k } \Big( \sum_{j\lesssim k}  2^{2s(j-k)} \, \Big)
   2^{2ks} \, \|u_k\|_{L^2}^2 \Big)
   \nonumber\\
   &\lesssim&
   \|D  u \|_{L^\infty} \, 
   \sum_j 2^{2js}   \|u_j\|_{L^2}^2  \,.
\eeqn
To pass to the second inequality, we used that
\eqn
   \|S_{j+1} Du \|_{L^\infty} \, 
   = \, \|m_{j+1}* Du\|_{L^\infty} 
   \, \lesssim \, \| Du \|_{L^\infty}  \, \|m_{j+1}\|_{L^1}
   \,,
\eeqn
where $\widehat{m_j}$ is the symbol of the Fourier multiplication operator $S_j$, and 
the fact that $\|m_j\|_{L^1}\sim1$ uniformly in $j$. 
Accordingly, we get
\eqn\label{eq-App-Enbd-2}
   \frac12 \partial_t \|u(t)\|_{\dot B^s_{2,2}}^2 \, \lesssim \, \|D  u(t)  \|_{L^\infty} \, 
   \|u(t)\|_{\dot B^s_{2,2}}^2 \,.
\eeqn
From
\eqn
    \|u(t)\|_{  B^s_{2,2}}^2 \, = \, \| u(t) \|_{L^2}^2 \, + \,
     \|u(t)\|_{\dot B^s_{2,2}}^2 \,,
\eeqn
and energy conservation, $\partial_t \| u(t) \|_{L^2}^2=0$, we 
obtain
\eqn\label{eq-App-Enbd-2}
   \frac12 \partial_t \|u(t)\|_{B^s_{2,2}}^2 \, 
   &=& 
   \frac12 \partial_t \|u(t)\|_{\dot B^s_{2,2}}^2 
   \nonumber\\ 
   &\lesssim& \|D  u(t)  \|_{L^\infty} \, 
   \|u(t)\|_{\dot B^s_{2,2}}^2 
   \nonumber\\
 & \lesssim & \|D  u(t)  \|_{L^\infty} \, 
   \|u(t)\|_{B^s_{2,2}}^2 \,.
\eeqn
This proves \eqref{eq-App-Enbd-1}. 
\qed

\subsection*{Acknowledgments}  
We are most grateful to Vlad Vicol for very useful comments. A
lso we would like to thank the referee for providing helpful remarks. 
The work of T.C. was supported by NSF grants DMS 1009448 and DMS 1151414 (CAREER). 
The work of N.P. was supported NSF grants DMS 0758247 and DMS 1101192 and an Alfred P. Sloan Research Fellowship.

\end{document}